\documentclass[12pt]{article}
\begin{document}
\thispagestyle{empty}
\begin{center}\Large{Permutations Containing and Avoiding 
$123$ and $132$ Patterns}
\end{center}
\vskip 10pt
\begin{center}
Aaron Robertson\footnote{webpage: 
www.math.temple.edu/\~{}aaron/\\ \indent 
This paper was supported in part by the NSF under the
PI-ship of Doron Zeilberger.} \\
Department of Mathematics\\ Temple University\\
Philadelphia, PA 19122\\
aaron@math.temple.edu
\end{center}
\vskip 10pt

\begin{abstract}
We prove that the number of permutations which avoid
$132$-patterns and have exactly one $123$-pattern, 
equals $(n-2)2^{n-3}$, for $n\geq 3$.
We then give a bijection onto the set of permutations 
which avoid $123$-patterns and have exactly one $132$-pattern.
Finally, we show that the number of permutations which
contain exactly one $123$-pattern and exactly one $132$-pattern
is $(n-3)(n-4)2^{n-5}$, for $n \geq 5$.
\end{abstract}

\begin{center}{\bf Introduction}\end{center}

In 1990, Herb Wilf asked the following:  How many permutations of
length $n$ avoid a given pattern, $p$?  By pattern-avoiding we
mean the following:  Let
$\pi$ be a permutation of length $n$ and let
$p=(p_1,p_2,\dots,p_k)$ be a permutation of length $k\leq n$
(we will call this a pattern of length $k$).  
Let $J$ be a set of $r$ integers, and let $j\in J$.
Define $place(j,J)$ to
be $1$ if $j$ is the smallest element in $J$, $2$ if it is the
second smallest, ..., and $r$ if it is the largest.
The permutation $\pi$ avoids the pattern $p$ if and only if
there does not exist a set of indices $I=(i_1,i_2,\dots,i_k)$, such that
$p=(place(\pi(i_1),I), place(\pi(i_2),I),\dots,place(\pi(i_k),I))$.

In two beautiful papers ([B] and [N]), the number of subsequences
containing exactly one $132$-pattern and exactly one $123$-pattern
are enumerated.  Noonan shows in [N] that the number of permutations
containing exactly one $123$-pattern is the simple formula
$\frac{3}{n}{{2n} \choose {n+3}}$.  B\'{o}na proves that the even
simpler formula ${{2n-3} \choose {n-3}}$ enumerates the number of
permutations containing exactly one $132$-pattern.  B\'{o}na's result 
proved a conjecture first made by Noonan and Zeilberger in [NZ].

Noonan and Zeilberger considered in [NZ] the number of permutations
of length $n$ which contain exactly $r$ p-patterns, for $r\geq 1$.
In this article we work towards the following generalization: 
How many permutations of length $n$ avoid patterns $p_i$, for
$i \geq 0$, and contain $r_j$ $p_j$-patterns, for $j \geq 1$, $r_j \geq 1$?
We will first consider the
permutations of length $n$ which avoid $132$-patterns, but contain
exactly one $123$-pattern.  We then define a natural bijection between
these permutations and the permutations of length $n$ which avoid $123$-patterns, 
but contain exactly one $132$-pattern.
Finally, we will calculate the number of permutations which contain
one $123$-pattern and one $132$-pattern.  These results address
questions first raised in [NZ].

\begin{center}{\bf Known Results}\end{center}
For completeness, two results which are already known are given below.

\noindent{\bf Lemma 1:  }{\it The number of permutations of length 
$n$ with one $12$-pattern is $n-1$.}

\noindent{\underline{Proof:}} Induct on $n$. The base case is trivial.
A permutation, $\phi$, of length $n$ with one $12$-pattern must have
$n=\phi(1)$ or $n=\phi(2)$.  If $n=\phi(1)$, by induction we get $n-2$
permutation.  If $n=\phi(2)$, then we must have $n-1=\phi(1)$ (or we
would have more than one $12$-pattern).  The rest of the entries of
$\phi$ must be decreasing.  Hence we get $1$ more permutation from
this second case, for a total of $n-1$.

\noindent{\bf Lemma 2:  }{\it The number of permutations which avoid
both the pattern $123$ and $132$ is $2^{n-1}$.}

\noindent{\underline{Proof:}} Let $f_n$ denote the number of permutations
we are interested in.  Then $f_n = \sum_{i=1}^{n} \, f_{n-i} + 1$ with $f_0=0$.
To see this, let $\rho$ be a permutation 
of length $n-1$.  Insert the element $n$ into the $i^{th}$ position of
$\rho$.  Call this new permutation of length $n$ $\pi$.  To assure that
$\pi$ avoids the $132$-pattern, we must have all entries preceding $n$ in
$\pi$ be larger that the entries following $n$.  To assure that
$\pi$ avoids the $123$-pattern, the entries preceding $n$ must be in
decreasing order.  This argument gives the sum in the recursion.  The
recursion holds by noting that if $n=1$, there is one permutation which
avoids both patterns.  To complete the proof note that $f_n=2^{n-1}$.

\begin{center}{\bf One 123-pattern, but no 132-pattern}\end{center}

\noindent{\bf Theorem 1:  }{\it The number of permutations of length $n$ which
have exactly one $123$-pattern, and avoid the $132$-pattern is
$(n-2)2^{n-3}$.}

\noindent{\underline{Proof:}} Let
$g_n$ denote the number of permutation we desire to count.  Call a
permutation {\it good} if it has exactly one $123$-pattern and avoids the
$132$-pattern.  Let $\gamma$ be a permutation of length $n-1$.
Insert the element $n$ into the $i^{th}$ position of $\gamma$. Call this
newly constructed permutation of length $n$, $\pi$.  To assure that
$\pi$ avoids the $132$ pattern, we must have all elements preceding $n$ 
in $\pi$ be larger than the elements following $n$ in $\pi$.  For
$\pi$ to be a {\it good} permutation, we must consider two
disjoint cases.

\noindent{\bf Case I:  } The pattern $123$ appears in the elements
following $n$ in $\pi$.  This forces the elements preceding $n$ to be
in decreasing order.  Summing over $i$, this case accounts for
$\sum_{i=1}^{n} \, g_{n-i}$ permutations.

\noindent{\bf Case II:  } The pattern $123$ appears in the elements
preceding and including $n$ in $\pi$.  
This forces the $3$ in the pattern to be $n$.  Hence the elements
preceding $n$ must contain exactly one $12$-pattern.
(Further there must be at least $2$ elements.  Hence
$i$ must be at least $3$).  From Lemma 1, 
this number is $i-2$.  We are also forced to avoid both
patterns in the elements following $n$.  Lemma 2 implies that
there are $2^{n-i-1}$ such permutations.  Summing over $i$, this
case accounts for $\sum_{i=3}^{n-1} \, (i-2)2^{n-i-1} \, +n-2$ permutations.

We have established that the recurrence relation

$$
g_n = \sum_{i=1}^{n} \, g_{n-i} + \sum_{i=3}^{n-1} \, (i-2)2^{n-i-1} \, +n-2, 
$$

\noindent
which holds for $n \geq 3$ ($g_0=0, g_1=0, g_2=0$), enumerates the
pemutations of length $n$ which avoid the pattern $132$ and contain
one $123$-pattern.

The obvious way to procede would be to find the generating function of
$g_n$.  However, in this article we would like to employ a
different, and in many circumstances more powerful, tool.  We will use
the Maple procedure {\tt findrec} in Doron Zeilberger's Maple
package {\tt EKHAD}\footnote{Available for download at 
{\tt www.math.temple.edu/\~{}zeilberg/}}.  (The Maple shareware
package {\tt gfun} could have also been used.)  Instructions for its use
are available online.  To use {\tt findrec} we compute the first few terms
of $g_n$.  These are (for $n\geq 4$) $4,12,32,80,192,448,1024$.  We
type {\tt findrec([4,12,32,80,192,448,1024],0,2,n,N)} and are given the
recurrence $h_n=4(h_{n-1} - h_{n-2})$ for $n \geq 4$.
Define $h_0=0, h_1=0, h_2=0,$ and $h_3=1$, and it is routine to verify that
$g_n=h_n$ for $n \geq 0$.  Another routine calculation shows us that
$h_n=(n-2)2^{n-3}$ for $n \geq 3$, thereby proving the statement of the 
theorem.

\begin{center}{\bf One 132-pattern, but no 123-pattern}\end{center}

\noindent{\bf Theorem 2:  }{\it The number of permutations of length $n$ which
have exactly one $132$-pattern, and avoid the $123$-pattern is
$(n-2)2^{n-3}$.}

\noindent
{\underline{Proof:}}
We prove this by exhibiting a (natural) bijection from the permutations
counted in Theorem 1 to the permutations counted in this theorem.
Define $S:=\{\pi : \pi$ avoids $132$-pattern and contains one $123$-pattern$\}$ and
$T:=\{\pi : \pi$ avoids $123$-pattern and contains one $132$-pattern$\}$.
We will show that $\mid S \mid = \mid T \mid$, by using the following 
bijection:

Let $\phi: S \longrightarrow T$.  Let $s \in S$, and 
let $abc$ be the $123$-pattern in $s$.  Then $\phi$ acts on the elements
of $s$ as follows:  $\phi(x) = x$ if $x \not \in \{b,c\}$, $\phi(b)=c$, and
$\phi(c)=b$.  In other words, all elements keep their positions except
$b$ and $c$ switch places.  An easy examination of several cases
shows that this is a bijection, thereby proving the theorem.

\begin{center}{\bf One 132-pattern and one 123-pattern}\end{center}

\noindent{\bf Theorem 3:  }{\it The number of permutations of length $n$ which
have exactly one $132$-pattern and one $123$-pattern is
$(n-3)(n-4)2^{n-5}$.}

\noindent
{\underline{Proof:}}
We use the same insertion technique as in the proof of Theorem 1.  
Let
$g_n$ denote the number of permutation we desire to count.  Call a
permutation {\it good} if it has exactly one $123$-pattern and exactly one
$132$-pattern.  Let $\gamma$ be a permutation of length $n-1$.
Insert the element $n$ into the $i^{th}$ position of $\gamma$. Call this
newly constructed permutation of length $n$, $\pi$. We note that the $132$-pattern
cannot consist of elements only preceding $n$.  If this were the case, we
would have two $123$-patterns ending with $n$. For
$\pi$ to be a {\it good} permutation, we must consider the following
disjoint cases.

\noindent{\bf Case I:  }  The $132$-pattern consists of elements following $n$.
In this case all elements preceding $n$ must be larger than the elements 
following $n$.
\vskip 5pt
\noindent
{\it Subcase A:  } The $123$-pattern consists of elements following $n$.  Summing
over $i$ we get $\sum_{i=1}^{n} \, g_{n-i}$ {\it good} permutations in this subcase.
\vskip 5pt
\noindent
{\it Subcase B:  } The elements preceding $n$ have exactly one $12$-pattern.
This gives a $123$-pattern where the 3 in the pattern is $n$.  We must also
avoid the $123$-pattern in the elements following $n$.  Summing over $i$ and
using Lemma 1 and Theorem 1, we get
$\sum_{i=3}^{n-3} \, (i-2)(n-i-3)2^{n-i-2}$ {\it good} permutations in this
subcase.

\noindent{\bf Case II:  } The $132$-pattern has the first element preceding
$n$, the last element following $n$, and $n$ as the middle element.
The elements preceding $n$ must be $n-1,n-2,\dots,n-1+2,n-i$, where
$n-i$ immediately precedes $n$ in $\pi$.
See
[B] for a more detailed argument as to why this must be true.
\vskip 5pt
\noindent
{\it Subcase A:  }  The elements preceding $n$ have exactly one $12$-pattern.
This gives a $123$-pattern where the last element of the pattern is $n$.
We must also avoid both the $123$ and the $132$ pattern in the elements following
$n$.  Summing over $i$ and using Lemma 1 and Lemma 2 we have
$\sum_{i=4}^{n-1} \, (i-3)2^{n-i-1}$ {\it good} permutations in this subcase.
\vskip 5pt
\noindent
{\it Subcase B:  }  The $123$-pattern consists of elements following $n$.  We
must have the elements preceding $n$ in $\pi$ be decreasing to avoid another
$123$-pattern.  Further, the elements following $n$ must not contain a $132$-pattern.
Using Theorem 1 and summing over $i$, we get a total of
$\sum_{i=2}^{n-3} \, (n-i-2)2^{n-i-3}$ {\it good} permutations in this
subcase.

In total, we find that the following recurrence enumerates the permutations
of length $n$ which contain exactly one $123$-pattern and one $132$-pattern.

$$
g_n = \sum_{i=1}^{n} \, g_{n-i} + \sum_{i=1}^{n-4} \, (2i(n-i-4)+n-3)2^{n-i-4}
$$

\noindent
for $n \geq 5$ and $g_1=g_2=g_3=g_4=0$.

Using {\tt findrec} 
again by typing
{\tt findrec([2,12,48,160,480,1344,3584],1,1,n,N)} \\ (where the list
is the first few terms of our recurrence for $n \geq 5$) we get
the recurrence $f_{n+1} = \frac{2(n+2)}{n} f_n$, with $f_1=2$.  After
reindexing, another routine calculation shows that $f_n=g_n$.  Solving
$f_n$ for an explicit answer, we find that $g_n=(n-3)(n-4)2^{n-5}$.

We conjecture that the number of 132-avoiding permutations with $r$
$123$-patterns is always a sum of powers of $2$.  For more evidence, and
further extensions see [ERZ].

\noindent
{\bf References}
\vskip 5pt
\noindent
[B] M. B\'{o}na, {\it Permutations with one or two $132$-subsequences.},
Discrete Mathematics, {\bf 181}, 1998, 267-274.
\vskip 2pt
\noindent
[ERZ] S. Ekhad, A. Robertson, D. Zeilberger, 
{\it The Number of Permutations With a Prescribed Number of
$132$ and $123$ Patterns}, submitted.  For a preprint see\\
{\tt www.math.temple.edu/\~{}[aaron,ekhad,zeilberg]/}.
\vskip 2pt
\noindent
[N] J. Noonan, {\it The Number of Permutations Containing Exactly One
Increasing Subsequence of Length $3$}, Discrete Mathematics, {\bf 152}, 
1996, 307-313.
\vskip 2pt
\noindent
[NZ] J. Noonan and D. Zeilberger, {\it The Enumeration of Permutations with
a Prescribed Number of ``Forbidden" Patterns}, Advances in Applied Mathematics,
{\bf 17}, 1996, 381-407.
\end{document}